\documentclass[a4paper,11pt]{article}
\usepackage{amsfonts,amsmath,amssymb,amstext,amsbsy,amsopn,amsthm}
\usepackage{url}
\usepackage{graphicx}

\newcommand{\D}{\displaystyle}
\newcommand{\itg}{\int \limits}
\newcommand{\R}{{\mathbf R}}
\newcommand{\N}{{\mathbf N}}
\newcommand{\Z}{{\mathbf Z}}
\newcommand{\cD}{\mathcal D}        
\newcommand{\cF}{\mathcal F}        
\newcommand{\cL}{\mathcal L}        
\newcommand{\cM}{\mathcal M}        
\newcommand{\cP}{\mathcal P}        
\newcommand{\cS}{\mathcal S}        
\newcommand{\eps}{\varepsilon}

\newcommand{\bdx}{\mathbf{x}}
\newcommand{\bdy}{\mathbf{y}}
\newcommand{\bdm}{\mathbf{m}}
\newcommand{\bdk}{\mathbf{k}}

\DeclareMathOperator{\Ree}{Re} 
\DeclareMathOperator{\Imm}{Im} 
\DeclareMathOperator{\erf}{erf} 
\DeclareMathOperator{\erfc}{erfc} 
\DeclareMathOperator{\wofz}{w} 
\DeclareMathOperator{\e}{e} 
 
\def\R{{\mathbb R}}
\def\Z{{\mathbb Z}}
\def\C{{\mathbb C}}

\newcommand{\keywords}{{\bf Keywords.  }}
\newcommand{\subjclass}{{\bf  Mathematics Subject Classification (2000). }}

\begin{document}
\title{On the fast computation of\\ high dimensional volume potentials}

\author{{\small Flavia Lanzara$^{\mbox{\tiny 1}}$ , Vladimir Maz'ya$^{\mbox{\tiny 2}}$,
Gunther Schmidt$^{\mbox{\tiny 3}}$}}
\date{}
\maketitle
\baselineskip=0.9\normalbaselineskip
\vspace{-7pt}
\hspace*{-5mm}
\parbox{10cm}{
\begin{flushleft}
{\footnotesize\em
\begin{itemize}
\item[$^{\mbox{\tiny\rm 1}}$]\rm Dipartimento di Matematica, Universit\`a
``La Sapienza'',\\
Piazzale Aldo Moro 2, 00185 Rome, Italy\\
\texttt{\rm lanzara\symbol{'100}mat.uniroma1.it}
\item[$^{\mbox{\tiny\rm 2}}$]Department of Mathematical Sciences, M\&O Building
University of Liverpool, Liverpool L69 3BX, UK\\Department of Mathematics, University of
Link\"oping, \\ 581 83 Link\"oping, Sweden\\
\texttt{\rm vlmaz\symbol{'100}mai.liu.se }
\item[$^{\mbox{\tiny\rm 3}}$]Weierstrass Institute for Applied Analysis and
Stochastics, \\  Mohrenstr. 39,
10117 Berlin, Germany \\
\texttt{\rm schmidt\symbol{'100}wias-berlin.de}
\end{itemize}
}
\end{flushleft}
}

\subjclass{65 D32, 65-05.}

 \keywords{Cubature of integral operators,
 multivariate approximation,  \\ separated representation }

\baselineskip=0.9\normalbaselineskip
\vspace{-3pt}

\begin{abstract}
A fast method of an arbitrary high order 
for approximating volume potentials is
proposed, which is effective
also in high dimensional cases. 
Basis functions
introduced
in the theory of approximate approximations are used. 
Results of numerical experiments,
which  show
approximation order $O(h^8)$ for the Newton potential 
in high dimensions, for example,
for $n= 200\,000$, are provided.
The computation time scales linearly in the space dimension.
New
one-dimensional integral representations
with separable integrands of
the potentials of advection-diffusion and heat equations are
obtained.
\end{abstract}

\section{Introduction}
\setcounter{equation}{0}
The cubature of high-dimensional volume potentials plays an important
role in a wide range of applications in physics, chemistry, biology,
financial mathematics etc.
Even a few years ago
this problem encountered unsurmountable difficulties
due to the so-called ``curse of dimensionality''.
With the development of separated representations (also called tensor-structured
approximations) by Beylkin and Mohlenkamp in  \cite{BM1,BM2}
the problem became tractable.
In fact, in recent years several 
fast algorithms for the computation of multi-dimensional convolutions
with singular kernels 
have  been proposed, cf. \cite{BCFH,GHB,H,HK,Khor,BHK}.   
They are mainly based on 
piecewise polynomial approximations of a separated
representation of the density.
Then a suitable separated approximation of the action
of the convolution operator on the basis functions allows
one to
determine the multi-dimensional convolution in question
by computing  a number of one-dimensional convolutions. 
On this way, the complexity of computing the multi-dimensional
convolutions on a  uniform tensor-product grid of size $h$
can be  reduced from $O(h^{-n} |\log h|)$, achieved with 
the traditional multi-variate FFT, to  $O(n h^{-1} |\log h|)$,
where $n$ is the space dimension.
A certain drawback of the schemes used in present
is the necessity
to find 
accurate separated representations
of the convolution operator acting on piecewise polynomials.
This procedure has been developed only for a
few particular kernels (\cite{BCFH,H}) and is rather involved, especially for higher order
approximations.

In this article we
propose 
a fast method of an arbitrary high order 
for approximating volume potentials, which is effective
also in high dimensional cases.
We use
basis functions
introduced
in the theory of {\it approximate approximations} \cite{Maz,MS2}, see also 
\cite {MSbook} and the references therein.
We report on  numerical experiments which  show
approximation order $O(h^8)$ for the Newton potential 
up to dimension $n= 200\,000$.
The computational complexity of the algorithm scales linearly
in the physical dimension. 
In the experiments,  high order cubatures 
lead to a minor
increase of the computation time
in comparison with second order approximations,
but are some orders of magnitude more accurate, especially in high dimension.

The method combines our basis functions with an approach 
from  Khoromskij \cite{Khor} for computing
volume potentials
on uniform or composite
refined grids.
For brevity of exposition we restrict ourselves to the  case
of uniform grids with step size $h$.
In \cite{Khor}, the density is approximated by piecewise constants and
the action of the convolution operator on the basis function
is written as a one-dimensional integral with a separable integrand, i.e.,
a product of functions depending only on one of the space
variables.
Then an accurate quadrature rule of this one-dimensional
integral provides a separated representation of the convolution
operator. The convergence orders $O(h^2)$ and
$O(h^3)$ using Richardson extrapolation are confirmed 
by numerical experiments for $n=3$.

Instead of piecewise constants we use
basis functions, which are Gaussians or products of Gaussians and
special polynomials and give rise to
high-order semi-analytic cubature formulas for
volume potentials.
In Section \ref{sec2} we describe these formulas 
and give error estimates.
According to \cite{MS5} (see also \cite[Section 6.3]{MSbook}), 
the action of volume potentials on the basis functions
allows for one-dimensional integral representations 
with separable integrands. 
In Section \ref{sec3} we demonstrate this by the example
of Newton's potential. We also obtain new one-dimensional  representations for 
the potentials of 
the advection-diffusion and the heat equations.
These one-dimensional integrals
in combination with a quadrature rule lead
to 
accurate separated representations
of the potentials.
In  Section \ref{numres},
we provide results of  numerical experiments
for the Newton potential showing that
even for very high space dimensions
these approximations preserve the
predicted convergence order of the cubature. 
In the final section \ref{quad}
we describe the quadrature of the one-dimensional integral
representations
used in the numerical examples.
Here we also
report
on
tests for second and fourth order formulas
of the Newton potential 
and the inverse of $-\Delta+a^2$,
which provide estimates of the rank of the separated representation
required to
approximate the action of the potential on a basis function
with a prescribed relative error.  

\section{Semi-analytic cubature formulas for  potentials}
\label{sec2}
\setcounter{equation}{0}
Here we collect some results
on high-order cubature formulas
for the volume potentials of the differential operators $-\Delta$ 
and $-\Delta+2 \mathbf{b} \cdot \nabla
+c$ and of the heat potential in $\R^n$.
More details can be found in \cite{MSbook}.

\subsection{Newton potentials in $\R^n$}
The Newton potential is
the inverse of the Laplace
operator and given in $\R^n$ by the formula
    \begin{equation}\label{harpot}
   {\cL}_n u(\bdx) =  \frac{\Gamma(\frac{n}{2}-1)}
   {4 \pi^{n/2}} \itg_{\R^n} \frac{u(\bdy)}{|\bdx-\bdy|^{n-2}} 
\, d\bdy \; , \quad  n \ge 3 \; .
    \end{equation}   
Here and in the following we denote by bold  $\bdx$ a vector of length
$n$ with the components $\bdx=(x_1,\ldots,x_n)$.
The integral \eqref{harpot} is a  unique solution of Poisson's equation
    \begin{equation*}
     -\Delta f = u \quad \textrm{in} \; \R^n \, , \quad |f(\bdx)| \le C
      |\bdx|^{n-2} \,  \;  \mbox{as} \; |\bdx| \rightarrow \infty \, .
    \end{equation*}   
In \cite{Maz,MS2} cubature formulas
for computing \eqref{harpot}
were constructed
which are based on the  approximation
of the density $u$ by
functions with analytically known Newton potentials.
In particular, one can approximate $u$ by the quasi-interpolant
    \begin{equation} \label{ch9eq2}
        u_h(\bdx) = {\cD}^{-n/2} \sum_{\bdm \in
{\Z}^n} u(h\bdm) \, \eta_{2M} \Big( \frac{\bdx-h\bdm}{\sqrt{\cD}h} 
\Big) \, ,
    \end{equation}
with the basis function
    \begin{equation}\label{LagM}
\eta_{2M}(\bdx)=\pi^{-n/2} \, L_{M-1}^{(n/2)}(|\bdx|^2) \, \e^{-|\bdx|^2} \, ,
\quad M \in \N \, ,
    \end{equation}   
where $L_{k}^{(\gamma)}$ are the generalized Laguerre polynomials 
    \begin{equation*} 
L_{k}^{(\gamma)}(y)=\frac{\e^{\,y} y^{-\gamma}}{k!} \, \Big(
{\frac{d}{dy}}\Big)^{k} \!
\left(\e^{\,-y} y^{k+\gamma}\right), \quad \gamma > -{\rm 1} \, .
    \end{equation*}
The Newton potential of $\eta_{2M}$ is given by
    \begin{equation} \label{ch3eta}
{\cL}_n \eta_{2M}(|\bdx|) =
\frac{1}{4 \pi^{n/2} |\bdx|^{n-2}}
\, \gamma\Big(\frac{n}{2} -1, |\bdx|^2\Big)
+ \frac{\e^{-|\bdx|^2}}{\pi^{n/2}}  
   \sum_{j=0}^{M-2}  \frac{L_{j}^{(n/2-1)}(|\bdx|^2)}{4 (j+1)} 
     \end{equation}
with the  incomplete Gamma function
    \begin{equation*}
\gamma(a,x)= \itg_0^x \tau^{a-1}  \e^{\, -\tau } \, d\tau \, .
    \end{equation*}               
Then the sum 
    \begin{equation}\label{conapp}
\cL_n u_h (\bdx) = \frac{h^2}{\cD^{n/2-1}} \sum_{\bdm \in
{\Z}^n} u(h\bdm) \cL_n \eta_{2M} \Big( \frac{\bdx-h\bdm}{\sqrt{\cD}h} 
\Big)
    \end{equation}   
is a semi-analytic cubature formula for $\cL_n u$.
It has been shown, that for sufficiently smooth and compactly supported
functions
\eqref{conapp} approximates $\cL_n u$
 with the error 
    \begin{equation}\label{harmest}
O(h^{2M})+ O(\e^{-\cD \pi^2} h^2) \, .
    \end{equation}   
This
follows from
the general asymptotic error estimate $O(h^{2M})+\eps$
for the quasi-interpolant \eqref{ch9eq2} with sufficiently smooth
and decaying generating functions $\eta$ satisfying the moment condition
   \begin{equation}  \label{momcon} 
    \itg _{{\R} ^n} \eta(\bdx) \, d \bdx = 1
  \; , \;                 
    \itg _{{\R} ^n} 
     \bdx^{\boldsymbol{\alpha}} \eta(\bdx) \, d \bdx = 0 , \quad
    \forall \, \boldsymbol{\alpha} \; , \; 
   1  \le |\boldsymbol{\alpha}| <  2M \, ,
    \end{equation}
see \cite[Chapter 2]{MSbook}.
The saturation error $\eps$ does not converge to zero as $h \to 0$, however,
because of 
\begin{align}\label{satur}
 \eps = O\Big( \max_{\bdx} \Big| \sum_{\bdk \in {\Z}^n  \backslash \{ \boldsymbol{0} \}}
\cF \eta (\sqrt{\cD} \bdk) \e^{\, 2 \pi i (\bdk,\bdx)/h}\Big| \Big) ,
\end{align}
it can be made arbitrarily small if the parameter $\cD$ is sufficiently large.
Additionally, the Newton potential maps
the fast oscillating saturation term \eqref{satur}
into a 
function with norm of order $O(h^2 \eps)$, which establishes the 
estimate \eqref{harmest}.
Therefore, in numerical computations with $\cD \ge 3$ 
formula \eqref{conapp}
behaves like
a high order
cubature formula.

\subsection{Potentials of advection-diffusion operators in $\R^n$}
The fundamental solution
of the operator $-\Delta +2 \mathbf{b} \cdot \nabla
+c$
with a vector $\mathbf{b} \in \C^n$ and $c \in \C$
depends on the value of $c + |\mathbf{b}|^2$. Here we use
the notation 
\[
\langle\mathbf{y},\mathbf{z}\rangle = \sum_{j=1}^n y_j z_j
\quad  \mbox{and} \quad |\mathbf{y}|^2 = \langle\mathbf{y},\mathbf{y}\rangle 
\]
also for vectors $\mathbf{y},\mathbf{z}\in \C^n$.
If $c + |\mathbf{b}|^2 \neq 0$, then the fundamental solution can be given as
\begin{align*}
\kappa_\lambda(\bdx)=
\frac{\e^{\langle\mathbf{b},\mathbf{x}\rangle}}{(2 \pi)^{n/2}}
\Big( \frac{|\bdx|}{\lambda}\Big)^{1-n/2}K_{n/2-1} (\lambda|\bdx|) \, ,
\end{align*}
where $\lambda \in \C \setminus (-\infty,0]$ with
$\lambda^2=c + |\mathbf{b}|^2$
and $K_\nu$ is the modified Bessel function of the second kind, also
known as Macdonald function, \cite[9.6]{Abr}. 
If $c + |\mathbf{b}|^2 = 0$, then for $n\ge 3$
\begin{align*}
\kappa_0(\bdx)=\frac{\Gamma(\frac{n}{2}-1)}
   {4 \pi^{n/2}} \frac{\e^{\langle\mathbf{b},\mathbf{x}\rangle}}{|\bdx|^{n-2}}
\end{align*}

To derive a cubature formula for the volume potential with the kernel
$\kappa_\lambda$
we look for a solution of
\begin{equation} \label{advec}
- \Delta f +2 \mathbf{b} \cdot \nabla f+c f =\e^{\, -|\bdx|^2} , \; \bdx \in \R^n \, ,
\end{equation}
which is given as the
one-dimensional integral
\begin{align} \label{pot_advdif_exp}
f(\bdx) &=\frac{ \e^{\,-|\mathbf{x}|^2} \, \lambda^{n/2-1}}
{|2\mathbf{x}+\mathbf{b}|^{n/2-1}}
\itg_0^\infty K_{n/2-1} (\lambda r) \, 
I_{n/2-1}(2|\mathbf{x}+\mathbf{b}|r)\,r\, \e^{\, -r^2}\, dr \,  ,
\end{align}
where $I_\nu$ is the modified Bessel function of the first kind, 
see \cite[Section 5.2]{MSbook}.

Using known analytic expressions of $I_{n+1/2}$ and $K_{n+1/2}$ (cf. \cite{Abr})
it is possible to derive analytic formulas of \eqref{pot_advdif_exp}
for odd space dimensions $n$. 
In particular, if $n=3$, then 
\begin{align*}
f(\bdx)= \frac{\sqrt{\pi}}{4} \frac{\e^{\,-|\mathbf{x}|^2}}{|2\mathbf{x}+\mathbf{b}|}
\Big(\wofz\Big(\frac{i}{2}\big(\lambda-|2\mathbf{x}+\mathbf{b}|\big) \Big) - 
\wofz\Big(\frac{i}{2}\big(\lambda+|2\mathbf{x}+\mathbf{b}|\big) \Big)\Big) \, ,
\end{align*}
where $\wofz$ denotes the scaled complementary error function
    \begin{equation*}
\wofz(z)=  \e^{\, -z^2} \, \erfc(-iz)
 =\e^{\, -z^2}\Big(1+ \frac{2 i}{\sqrt{\pi}}  \itg_{0}^{z} \e^{t^2} \,
 dt \Big)  
    \end{equation*}
and  $\erfc = 1- \erf$ is the complementary error function.
From
    \begin{equation} \label{LaguerM} 
 L_{M-1}^{(n/2)}(|\mathbf{x}|^2) 
        \, \e^{-|\mathbf{x}|^2} 
= \sum_{k=0}^{M-1}
        \frac{(-1)^k}   {k! \, 4^k} \, \Delta^k  \e^{-|\mathbf{x}|^2}
    \end{equation}
(see \cite[Theorem 3.5]{MSbook}), one can derive as in the case of Newton potentials
semi-analytic cubature formulas for the potential of the
advection-diffusion equation
with the approximation rate \eqref{harmest}.

\subsection{Heat potential} \label{scnheat}
Next we consider the non-homogeneous
heat equation
\begin{equation}\label{heatn}
\begin{split}
 f_{t}-\nu   \Delta_{\bdx} f &=u(\bdx,t)\hskip1cm
\> t\geq 0,\> \nu > 0 \, ,\\
  f(\bdx,0)&=0 \hskip2cm \bdx\in \R^{n}.
\end{split}
\end{equation}
It is well known that 
the solution of this Cauchy problem
can be written as
\begin{equation*}
f(\bdx,t)=\itg_0^t (\cP_{t-\lambda} u(\cdot,\lambda))(\bdx) \, d\lambda \, ,
\end{equation*}
where $\cP_t$ is the Poisson integral
\begin{equation*}
(\cP_t   u(\cdot,\lambda))(\bdx)=\frac{1}{({4 \pi\nu t)^{n/2} }} \itg_{\R^n} 
    {\e}^{-|\bdx-\bdy|^2/(4\nu t)}  u(\bdy,\lambda) \, d\bdy \, . 
\end{equation*}
Suppose that the right-hand side $u$ is a sufficiently smooth
function with compact support in $\R^n \times \R_+$.
An approximation of the solution $f(\bdx,t)$ of \eqref{heatn} can be obtained if 
$u$ is approximated
 by the quasi-interpolant on the
rectangular grid $\{(h\bdm,\tau i)\}$
\begin{equation}\label{fht2}
    u_{h,\tau}(\bdx,t)= 
    \frac{\pi^{-(n+1)/2}}{\sqrt{\cD_{0}\cD^{n}}} 
\! \sum_{\substack{i\in\Z\\ \bdm\in\Z^{n}}}\!
  u(h\bdm,\tau i) \e^{-(t- \tau i)^{2}/(\cD_{0}\tau^{2})}
    \e^{-|\bdx-h\bdm|^{2}/(\cD h^{2})}.
\end{equation}
It can be easily seen, that
\begin{align*}
(\cP_{t-\lambda}u_{h,\tau}(\cdot,\lambda))(\bdx)
=\frac{h^{n}}{\pi^{n/2}\sqrt{\pi \cD_{0}}}
\sum_{\substack{i\in\Z\\ \bdm\in\Z^{n}}}
u(h\bdm,\tau i)\frac{\e^{-(\lambda- \tau i)^{2}/(\cD_{0}\tau^{2})} 
\e^{\, - |\bdx-h\bdm|^2/(\cD h^{2}+4\nu (t-\lambda))}}{(\cD h^{2} +
   4\nu  (t-\lambda) )^{n/2}},
\end{align*}
hence the sum
\begin{align*}
   f_{h,\tau}(\bdx,t)= \itg_0^t  (\cP_{t-\lambda}
   u_{h,\tau}(\cdot,\lambda))(\bdx) \, d\lambda
  = \frac{h^{n}}{\pi^{n/2}\sqrt{\pi \cD_{0}}}
    \sum_{\substack{i\in\Z\\ \bdm\in\Z^{n}}}
     K_2(\bdx-h\bdm,t,\tau i)\, u(h\bdm,\tau i)
\end{align*}
is an approximation of $f(\bdx,t)$, where we use the notation
\begin{equation} \label{heatint}
K_2(\bdx,t,\tau i)=\itg_{0}^{t}\frac{ \e^{-(t-\lambda-\tau i)^{2}/(\cD_{0}\tau^{2})}
   \e^{-|\bdx|^{2}/(\cD h^{2} +4\nu \lambda)}}{(\cD h^{2} +
   4\nu  \lambda )^{n/2}} \, d\lambda \, .
\end{equation}
Since  $u_{h,\tau}(\bdx,t)$ differs from $u(\bdx,t)$ by 
\begin{align*}
     \D | u(\bdx,t)-u_{h,\tau}(\bdx,t)|
   \D \leq \eps +c \big( (\tau \sqrt{\cD_{0}})^{2}+(h \sqrt{\cD})^{2}\big)\, , \quad \forall \bdx\in 
    \R^{n},\, t\in [0,T] \, ,
\end{align*}
the function $f_{h,\tau}(\bdx,t)$ approximates the solution
$f(\bdx,t)$ with the error 
\begin{align*}
    |f(\bdx,t)-f_{h,\tau}(\bdx,t)| &= \D \frac{1}{({4 \pi\nu)^{n/2} }}
\itg_{0}^{t} \, d\lambda
    \itg_{\R^{n}}\frac{\e^{-|\bdx-\bdy|^{2}/(4(t-\lambda))}}{(t-\lambda)^{n/2}}
   | u(\bdy,\lambda)-u_{h,\tau}(\bdy,\lambda)|\, d\bdy\\  
   & \leq T\, ||u-u_{h,\tau}||_{L^{\infty}(\R^{n} \times
  [0,T])} \leq \D \eps +c\big( (\tau \sqrt{\cD_{0}})^{2}+(h
\sqrt{\cD})^{2}\big)
\end{align*}
for all $\bdx\in \R^{n}$ and $t\in [0,T]$.

 \section{One-dimensional integral representations of the  
potentials acting on Gaussians}\label{sec3}
\setcounter{equation}{0}
The computation of the mentioned cubature formulas 
on the grid $\{h\bdk\}$
leads to discrete convolutions
    \begin{equation} \label{conapp1}
\sum _{\bdm}
a_{\bdk-\bdm} \, u(h\bdm)  \, ,
    \end{equation}   
where the indices $\bdm,\bdk$ belong to some subset of $\Z^n$, and the
coefficients $a_{\bdk}$ are given either analytically 
or by a one-dimensional integral with smooth 
integrand.
For general functions $u$  the most efficient summation
methods
for \eqref{conapp1} are probably
fast  convolutions
based on multi-variate FFT.
However, even for the space dimension $n=3$
problems of moderate size often exceed  the capacity of
available computer systems.

The situation in high dimensions is much better
when the vectors $\{u(h\bdm)\}$ and $\{a_{\bdk}\}$
allow {\it separated representations} \cite{BM1}, i.e. for given accuracy $\epsilon$
they can be represented as a sum of products of vectors in dimension $1$
\begin{align*} 
&u(hm_1,\ldots,hm_n) = \sum_{p =1}^R r_p \prod_{j=1}^{n}
u_j^{(p)}(hm_j) + O(\epsilon) \, , \\
&a(k_1,\ldots,k_n) = \sum_{p =1}^R s_p \prod_{j=1}^{n}
v_j^{(p)}(k_j) + O(\epsilon)\, .
\end{align*} 
Then an approximate value of \eqref{conapp1} can be computed by the sum of
products of one-dimensional convolutions
\begin{align*} 
\sum _{\bdm}  
a_{\bdk-\bdm} \, u(h\bdm) \approx \sum_{p,q=1}^R r_p s_q
\prod_{j=1}^{n}\sum_{m_j}  v_j^{(q)}(k_j-m_j)u_j^{(p)}(hm_j)\, .
\end{align*} 
To obtain a separated representation of the vector $\{a_{\bdk}\}$
the following idea from  Khoromskij \cite{Khor} can be applied: 
Suppose that the density $u$ is approximated by interpolation or
quasi-interpolation using linear combinations of the translates $\phi(\cdot -h\bdm)$ of a basis
function.
Then the components  of $\{a_{\bdk}\}$ are the values  of the
volume potential with the
kernel $g$
acting on $\phi$,
\[
a_{\bdk}= \itg_{\R^n}g(\bdk-\bdy ) \phi(\bdy)\, d \bdy \, .
\] 
If this integral can be transformed to
a one-dimensional
integral
with {\it separable integrand}
\[
\itg_{\R^n} g(\bdx-\bdy ) \phi(\bdy)\, d \bdy 
= \itg_0^\infty \prod_{j=1}^{n} g_j(x_j,t) \, dt \, ,
\] 
then a separated approximation of $\{a_{\bdk}\}$ is obtained by
applying an accurate quadrature rule for 
this integral 
\[
a_{\bdk}\approx \sum_{\ell=1}^R \omega_\ell\prod_{j=1}^{n}g_j(k_j,\tau_\ell)\, .
\]
In \cite{Khor} one-dimensional
integral representations
were obtained for  the potentials of the Laplace and of the  advection-diffusion equation
with $\mathbf{b}=0$, $c > 0$, applied to the characteristic
functions $\phi$ of a cube in $\R^3$.
This gives rise to a fast algorithm which is based on piecewise constant basis functions.
The convergence order $O(h^2)$  is  proved together with 
the order $O(h^3)$ by using Richardson extrapolation.

In this section we describe 
one-dimensional integral representations with separable integrands
for
the Newton potential and 
the potentials of 
the advection-diffusion and the heat equations
acting on basis functions 
introduced in the previous section.

\subsection{Newton potential}

The separated representation
of the  second order cubature formula 
    \begin{equation}\label{ch9harm3}  
   \frac{{\cD}h^2}{(\pi {\cD})^{n/2}}
   \sum_{\mathbf{m} \in{\Z}^n} u(h\mathbf{m}) \, \cL_n
   \big(\e^{-|\,\cdot\,|^2}\big)(\mathbf{r_m})
    \end{equation}                                     
for the Newton potential with
    \begin{equation*}
\mathbf{r_m} = \frac{\mathbf{x}-h\bdm}{\sqrt{\cD}h} 
\quad \mbox{and} \quad {\cL}_n \big(\e^{\,-|\,\cdot\,|^2}\big)(\bdx)
= \frac{1}{4 |\bdx|^{n-2}} \, \gamma\Big(\frac{n}{2} -1, {|\bdx|^2}\Big)\, ,
    \end{equation*}
follows from the formula 
    \begin{equation*}
{\cL}_n \big(\e^{\,-|\,\cdot\,|^2}\big)(\bdx)
= \frac{1}{4}
\itg_0^\infty \frac{\e^{\, - |\bdx|^2/(1+t)}}{(1+ t)^{n/2}}
\, dt  \, ,
    \end{equation*}
which was obtained in \cite{MS5} and is valid for $n \ge 3$.
The quadrature rule of the last integral
    \begin{equation*}
{\cL}_n \big(\e^{\,-|\,\cdot\,|^2}\big)(\bdx)
\approx \frac{1}{4}\sum_{p=1}^R \omega_p \frac{\e^{\, -
    |\bdx|^2/(1+\tau_p)}}{(1+ \tau_p)^{n/2}} 
= \frac{1}{4}\sum_{p=1}^R \omega_p\prod_{j=1}^{n}
\frac{\e^{\, -  x_j^2/(1+\tau_p)}}{(1+ \tau_p)^{1/2}}
    \end{equation*}
with certain quadrature  weights $\omega_p$ and nodes $\tau_p$ 
gives
the separated approximation
\[
\cL_n u_h (h \bdk) \approx\frac{\cD h^2}{(\pi\cD)^{n/2} }  \sum_{\bdm \in
{\Z}^n} u(h\bdm)
\sum_{k=1}^R \omega_k \prod_{j=1}^{n}
\frac{\e^{- (k_j- m_j)^2/(\cD (1+\tau_k))}}
        {(1+\tau_k)^{1/2}} \, .
\]

To get an
approximation for higher order cubature
formulas we note that  the same convergence order
$O(h^{2M})+ O(\e^{-\cD \pi^2} h^2)$ as in the case of generating functions \eqref{LagM}
holds,
when the density is approximated by
the sum 
\begin{align*}
\cM_M u(\bdx)= (\pi \cD)^{-n/2}
\sum_{\bdm \in \Z^n} u(h\bdm)
\prod_{j=1}^n \widetilde \eta_{2M} \Big( \frac{x_j - h m_j}{\sqrt{\cD} h}\Big) .
 \end{align*}
Here the basis function is the tensor product of
the one-dimensional  functions
\begin{align} \label{Lagone}
\widetilde \eta_{2M} (x)= L_{M-1}^{(1/2)}(x^2)
\e^{\, - x^2} \, ,
 \end{align}
which obviously satisfies the moment condition \eqref{momcon}. 
Thus the cubature formula
\begin{align} \label{cubfM}
\cL_{M,h}^{(n)}u(\bdx)=\frac{h^2}{\cD^{n/2-1}} \sum_{\bdm \in
{\Z}^n} u(h\bdm) \cL_n  \Big(\prod_{j=1}^n \widetilde \eta_{2M} \Big) \Big( \frac{\bdx-h\bdm}{\sqrt{\cD}h} 
\Big)
 \end{align}
approximates $\cL_n u$ with the error \eqref{harmest}. 
The relation
\begin{align} \label{1dsum}
\widetilde\eta_{2M} (x)&=
\sum_{k=0}^{M-1}
        \frac{(-1)^k}   {k! \, 4^k} \, 
\frac{d^{2k}}{dx^{2k}}  \e^{-x^2} 
=\e^{\, - x^2} \sum_{k=0}^{M-1}
L_{k}^{(-1/2)}(x^2)\, ,
 \end{align}
leads to the one-dimensional integral representation
of $\cL_n \Big(\prod_{j=1}^n \widetilde \eta_{2M}\Big)$
by writing the solution of the Poisson equation
\begin{align*}
- \Delta u(\bdx)
= \prod_{j=1}^n \widetilde \eta_{2M} (x_j)
 \end{align*}
as the integral
\begin{align*}
u(\bdx)&=\frac{1}{4} \, \prod_{j=1}^n \sum_{k=0}^{M-1}
         \frac{(-1)^k }   {k! \, 4^k} \, 
\frac{d^{2k}}{dx_j^{2k}}\itg_0^\infty 
\frac{\e^{- |\bdx|^2/(1+t)}}
        {(1+t)^{1/2}} \, dt  \\
&=\frac{1}{4 } \itg_0^\infty 
\Big(\prod_{j=1}^n \sum_{k=0}^{M-1}
         \frac{(-1)^k }   {k! \, 4^k} \, 
\frac{d^{2k}}{dx_j^{2k}}\e^{-x_j^2/(1+t)}\Big)
\frac{ dt } {(1+t)^{n/2}} \\
&=\frac{1}{4  } \itg_0^\infty 
\Big(\prod_{j=1}^n  \sum_{k=0}^{M-1}
        \frac{\e^{\, -x_j^2/(1+t)}}   {(1+t)^{k}} \, 
L_{k}^{(-1/2)}\Big(\frac{x^2_j}{1+t}\Big)
\Big) \,  \frac{ dt } {(1+t)^{n/2}} \, .
 \end{align*}

Hence the separated representation
of the Newton potential
is reduced to find accurate quadrature rules for the
parameter dependent integrals 
\begin{equation} \label{newton2}
I_M(\bdx)=\itg_0^\infty 
\Big(\prod_{j=1}^n \e^{\, -x_j^2/(1+t)}\sum_{k=0}^{M-1}
        \frac{1 }   {(1+t)^{k}} \, 
L_{k}^{(-1/2)}\Big(\frac{x^2_j}{1+t}\Big)
\Big) \,  \frac{ dt } {(1+t)^{n/2}}  \, , \; M \ge
1 \,.
\end{equation}
It is necessary 
that the quadrature rule
approximates the
integrals \eqref{newton2}
with prescribed error for the parameters $x_j= (k_j-m_j)/\sqrt{\cD}$ within the range
$|x_j| \le K$ and some given bound $K$.
This will be discussed in Section \ref{quad}.

\subsection{Potentials of advection-diffusion operators}
Here we look for a one-dimensional integral representation with
separable integrand
of the 
integral \eqref{pot_advdif_exp}.
Let $f$ be the solution of the advection-diffusion equation \eqref{advec}.
Then  $v = f \e^{-\langle\mathbf{b},\mathbf{x}\rangle}$
satisfies
\begin{align*}
-\Delta v +  \lambda^2 v = \e^{-|\bdx+\mathbf{b}/2|^2}  
\e^{|\mathbf{b}|^2/4}  \quad \mbox{with} \;\lambda^2 = c + |\mathbf{b}|^2
\end{align*}
in $\R^n$. In was shown in  \cite{MS5} that
$v$ can be represented in the form 
\begin{align} \label{1drepadv}
v(\bdx)=\frac{\e^{|\mathbf{b}|^2/4} }{4}
\itg_0^\infty \frac{\e^{-\lambda^2 t/4}
\e^{\, - |\bdx+\mathbf{b}/2|^2/(1+t)}}{(1+ t)^{n/2}}
\, dt 
\end{align}
for 
$\Ree \lambda^2 \ge 0$ and $n \ge 3$.
If $\Ree \lambda^2 > 0$, then \eqref{1drepadv} is valid for all
space dimensions $n$,
see also \cite[Theorem 6.4]{MSbook}.
Thus 
\begin{align}
f(\bdx)&=\itg_{\R^n}\kappa_\lambda(\bdx-\bdy) 
\e^{\, -|\bdy|^2} \, d \bdy
=\frac{\e^{ \langle\mathbf{b},\mathbf{x}\rangle}
\e^{|\mathbf{b}|^2/4 }}{4} 
\itg_0^\infty \frac{\e^{-t(c + |\mathbf{b}|^2) /4}
\e^{\, - |\bdx+\mathbf{b}/2|^2/(1+t)}}{(1+ t)^{n/2}}
\, dt \nonumber \\
&=\frac{1 }{2}
\itg_0^\infty \frac{\e^{-c \, t/2}}
{(1+ 2t)^{n/2}}\e^{\, - |\bdx-t \mathbf{b}|^2/(1+2t)}
\, dt \, . \label{intadv}
\end{align}
Consequently, if $\Ree (c + |\mathbf{b}|^2) \ge 0$, then an
approximate  solution of
the  advection-diffusion equation 
\[
- \Delta f +2 \mathbf{b} \cdot \nabla f+c f  = u
\]
is given in $\R^n$ by the sum
\begin{align*}
f_h(\bdx)=
\frac{\cD h^2 }{2 (\pi\cD)^{n/2} }
 \sum_{\bdm \in  {\Z}^n}u(h
 \bdm)
\itg_0^\infty 
\frac{\e^{-\cD h^2  c \, t/2}}
        {(1+2t)^{n/2}}  
\e^{- |\mathbf{x}-h\bdm- t\cD h^2 \mathbf{b}|^2/(\cD h^2(1+2t))}\, dt \, ,
 \end{align*}
which approximates  $f$ with the order $O(h^2)$.
Analogously to the case of Newton potentials 
the one-dimensional integrals  for cubature formulas of order
\eqref{harmest}  are based on the basis functions 
$\prod \widetilde \eta_{2M}(x_j)$.
Now the relation \eqref{1dsum} leads to
\begin{align}
\itg_{\R^n} \kappa_\lambda(\bdx- \bdy)&
\prod_{j=1}^n \widetilde \eta_{2M}(y_j) \, d \bdy
= \prod_{j=1}^n \sum_{k=0}^{M-1}
         \frac{(-1)^k }   {k! \, 4^k} \, 
\frac{\partial^{2k}}{\partial x_j^{2k}} \itg_{\R^n} 
\kappa_\lambda(\bdx-\bdy) \e^{- |\bdy|^2}\, d \bdy \nonumber \\
&= \frac{ 1}{2}  
\itg_0^\infty \frac{\e^{-c\,  t/2}}
        {(1+2t)^{n/2}} 
 \Big( \prod_{j=1}^n \sum_{k=0}^{M-1}
         \frac{(-1)^k }   
{k! \, 4^k}\frac{d^{2k}}{dx_j^{2k}}\e^{\, - |\bdx-t \mathbf{b}|^2/(1+2t)} \Big) dt \nonumber \\
&=\frac{ 1}{2  } \itg_0^\infty \frac{\e^{-c\,  t/2}}
        {(1+2t)^{n/2}}
\Big(\prod_{j=1}^n g_M(2t,x_j-t b_j) \Big) dt \label{intadvM}
 \end{align}
with the function
\[
g_M(t,x) =   \sum_{k=0}^{M-1}
        \frac{\e^{\, -x^2/(1+t)}}   {(1+t)^{k}} \, 
L_{k}^{(-1/2)}\Big(\frac{x^2}{1+t}\Big).
\]

\subsection{Heat potential}
Recall that the approximation of the non-homogeneous heat equation
\eqref{heatn} on the
rectangular grid $\{(h \bdk,\tau \ell)\}$ is given by the
$(n+1)$-dimensional convolution
\begin{align*}
   f_{h,\tau}(h \bdk, \tau \ell)
=  \frac{h^{n}}{\pi^{(n+1)/2}\sqrt{\cD_{0}}}
    \sum_{\substack{i\in\Z\\ \bdm\in\Z^{n}}}
     K_2(h(\bdk-\bdm),\tau\ell,\tau i)\, u(h\bdm,\tau i) \, ,
\end{align*}
where the integral $K_2$ given by 
\eqref{heatint} cannot be taken analytically.
Making the substitution
\[
\lambda= \frac{t}{1+\e^{-\xi}}
\]
$K_2(\bdx,t,\tau i)$ transforms to the integral over $\R$
\begin{align*}
K_2(\bdx,t,\tau i)&=\frac{t}{4 }
\itg_{-\infty}^{\infty}\frac{ \e^{-(t/(1+\e^{\xi})-\tau i)^{2}/(\cD_{0}\tau^{2})}
   \e^{-|\bdx|^{2}/(\cD h^{2} +4\nu t/(1+\e^{-\xi}))}}
{(\cD h^{2} +   4\nu  t/(1+\e^{-\xi}) )^{n/2}\cosh^2 (\xi/2)} \, d \xi 
\end{align*}
with exponentially decaying integrand.
The quadrature rule for $K_2(h\bdk,\tau \ell,\tau i)$
\[
\frac{\tau \ell}{4 }
\sum_{p=1}^{R}
\omega_p\, \frac{ \e^{-(\ell/(1+\e^{\xi_p}) -i )^{2}/\cD_{0}}}
{(\cD h^{2} +   4\nu  \tau \ell/(1+\e^{-\xi_p}) )^{n/2}\cosh^2 (\xi_p/2)}  
\prod_{j=1}^{n} \e^{-k_j^{2}/(\cD  +4\nu \tau \ell/(h^{2}(1+\e^{-\xi_p})))}
\]
with certain  weights $\omega_p$ and nodes $\xi_p$
 and the separated representation
of the right-hand side of \eqref{heatn}
\begin{align*} 
u(h\bdx,t) = \sum_{q =1}^Q r_q \prod_{j=1}^{n}
u_j^{(q)}(x_j,t) + O(\epsilon) 
\end{align*} 
leads to the approximation
\begin{align*} 
   f_{h,\tau}(h \bdk, \tau \ell) \approx 
 \frac{h^{n}\tau \ell}{4\pi^{(n+1)/2}\sqrt{\cD_{0}}} 
 \sum_{q =1}^Q \sum_{p=1}^{R}  \frac{\omega_p\, r_q}{\cosh^2 (\xi_p/2)}  \sum_{i\in\Z}
a_{\ell,i,p }\, c^{(q)}_{\bdk,\ell,i,p} \, ,
 \end{align*}
where we use the abbreviation
\begin{align*} 
&a_{\ell,i,p}=
 \frac{ \e^{-(\ell/(1+\e^{\xi_p}) -i )^{2}/\cD_{0}}}
{(\cD h^{2} +   4\nu  \tau \ell/(1+\e^{-\xi_p}) )^{n/2}} \, ,
 \end{align*} 
and $c^{(q)}_{\bdk,\ell,i,p}$ is the product of one-dimensional convolutions
 \begin{align*} 
&c^{(q)}_{\bdk,\ell,i,p}= \prod_{j=1}^{n}
\sum_{m_j \in \Z} u_j^{(q)}(hm_j,\tau i)
\e^{-(k_j-m_j)^{2}/(\cD  +4\nu \tau \ell/(h^{2}(1+\e^{-\xi_p})))}\, .
\end{align*} 

Approximations which converge with higher order to the solution of \eqref{heatn}
can be obtained  if the right-hand side $u$ is approximated by 
\begin{align} \label{fht}
u_{h,\tau}(\bdx,t) =   \frac{\pi^{-(n+1)/2}}{\sqrt{\cD_{0}\cD^{n}}}
\sum_{\substack{i\in\Z\\ \bdm\in\Z^{n}}}
  u(h\bdm,\tau i)
\widetilde \eta_{2S}\Big( \frac{t - \tau i}{\sqrt{\cD_0}  \tau}\Big) 
\prod_{j=1}^n \widetilde \eta_{2M} \Big( \frac{x_j - h m_j}{\sqrt{\cD} h}\Big)
 \end{align}
with the basis functions $\widetilde \eta_{2M}$  defined by \eqref{Lagone}.
Then by using \eqref{1dsum}
the heat potential of the quasi-interpolant \eqref{fht} is given by
\begin{align*}
   f_{h,\tau}(\bdx,t)= \itg_0^t  (\cP_{t-\lambda}
   u_{h,\tau}(\cdot,\lambda))(\bdx) \, d\lambda  
  = \frac{h^{n}}{\pi^{n/2}\sqrt{\pi \cD_{0}}}
    \sum_{\substack{i\in\Z\\ \bdm\in\Z^{n}}}
     K_{S,M}(\bdx-h\bdm,t,\tau i)\, u(h\bdm,\tau i) \, , 
\end{align*}
where we denote
 \begin{align*}
K_{S,M}(\bdx,t,\tau i)=\sum_{k=0}^{S-1}
         \frac{(-1)^k(\cD_0 \tau^2)^k }   {k! \, 4^k} \, \frac{\partial^{2k}}{\partial t^{2k}} 
\prod_{j=1}^{n} \sum_{\ell=0}^{M-1}
         \frac{(-1)^\ell(\cD h^2)^\ell }   {\ell! \, 4^\ell} \,
\frac{\partial^{2\ell}}{\partial x_j^{2\ell}}
K_2(\bdx,t,\tau i)\, .
\end{align*}
We have
 \begin{align*}
\frac{(-1)^\ell}{\ell! \, 4^\ell}
\frac{\partial^{2\ell}}{\partial x_j^{2\ell}} 
\e^{-|\bdx-h\bdm|^{2}/(\cD h^{2} +4\nu \lambda)}=  
\frac{\e^{-|\bdx-h\bdm|^{2}/(\cD h^{2} +4\nu \lambda)}}{(\cD h^{2} +4\nu \lambda)^\ell} 
L_{\ell}^{(-1/2)}\Big(\frac{(x_j-h\,m_j)^2}{\cD h^{2} +4\nu \lambda}\Big) \, ,
\end{align*}
which by using \eqref{heatint} leads to the representation
with separable integrand
 \begin{align*}
K_{S,M}(\bdx,t,\tau i) &=
\sum_{k=0}^{S-1}
         \frac{(-1)^k(\cD_0 \tau^2)^k }   {k! \, 4^k} \, 
\frac{\partial^{2k}}{\partial t^{2k}} 
\itg_{0}^{t}{ \e^{-(\lambda-(t- \tau i))^{2}/(\cD_{0}\tau^{2})}} 
\prod_{j=1}^{n} g_M(\lambda,x_j) d\lambda \, .
\end{align*}
Here we denote by $g_M$ the function
 \begin{align*}
g_M(\lambda,x)=    \e^{-x^{2}/(\cD h^{2} +4\nu \lambda)}\sum_{\ell=0}^{M-1}
\frac{(\cD h^2)^\ell}{(\cD h^{2} +4\nu \lambda)^{\ell+1/2}}
L_{\ell}^{(-1/2)}\Big(\frac{x^2}{\cD h^{2} +4\nu \lambda}\Big) 
\, .
\end{align*}

\section{Numerical results for the Newton potential} \label{numres}
\setcounter{table}{0}
We provide results of some experiments 
which show the accuracy and numerical
convergence orders of the method.
In the cube $[-6,6]^n$, $n \ge 3$, we compute the Newton potential
of the densities
\[
u_1(\bdx) = \e^{- |\bdx|^2} \quad \mbox{and} \quad
u_2(\bdx) = (4|\bdx|^2 -2n)\e^{- |\bdx|^2} \, ,
\] 
which have the exact values
\[
{\cL}_nu_1(\bdx) = \frac{1}{4 |\bdx|^{n-2}} \, 
\gamma\Big(\frac{n}{2} -1,|\bdx|^2\Big)
\quad \mbox{and} \quad
{\cL}_n u_2(\bdx) = - \e^{- |\bdx|^2} \, .
\]

In Table~\ref{tabel1}
we compare the values of the exact and  
the approximate solution for ${\cL}_n u_1$
at some points $(x_1,0,\ldots,0) \in \R^n$
of the grid, where we have chosen the space dimensions
$n=3,10,100,300$.

\begin{table}[!h]
\begin{small}
\begin{center}
\begin{tabular}{c|c|cc|c}
dimension&  $\hskip3mm x_1$ \hspace{3mm} & \hspace{12mm}  exact \hspace{12mm} & 
\hspace{5mm} approximation \hspace{5mm} & 
\hskip2mm relative error \hspace{2mm} \\ \hline 
 &   0.0 &    0.50000000000000 &    0.49999999923850 &  1.5230E-09 \\
 &   1.0 &    0.37341206640621 &    0.37341206614375 &  7.0287E-10 \\
$n=3$ &   2.0 &    0.22052034769061 &    0.22052034766043 &  1.3685E-10 \\
 &   3.0 &    0.14770122470992 &    0.14770122470423 &  3.8549E-11 \\
 &   4.0 &    0.11077836397370 &    0.11077836396658 &  6.4242E-11 \\
 &   5.0 &    0.08862269254514 &    0.08862269253834 &  7.6764E-11 \\
\hline \hline 
&   0.0 &    0.06250000000000 &    0.06249999932963 &  1.0726E-08  \\
 &   1.0 &    0.02848223531423 &    0.02848223504590 &  9.4209E-09 \\
$n=10$ &   2.0 &    0.00331951101348 &    0.00331951099712 &  4.9280E-09 \\
 &   3.0 &    0.00022377080789 &    0.00022377080851 &  2.7741E-09 \\
 &   4.0 &    0.00002288605175 &    0.00002288605169 &  2.6127E-09 \\
 &   5.0 &    0.00000383999984 &    0.00000383999985 &  6.8146E-10 \\
  \hline \hline 
 &  0.0 &    0.00510204081633 &    0.00510204386664 &  5.9786E-07 \\
 &   1.0 &    0.00191522512312 &    0.00191522620272 &  5.6369E-07 \\
$ n=100$ &   2.0 &    0.00010155802170 &    0.00010155808096 &  5.8347E-07 \\
 &   3.0 &    0.00000076714427 &    0.00000076714503 &  9.9929E-07 \\
 &   4.0 &    0.00000000084085 &    0.00000000084085 &  1.8801E-06 \\
 &   5.0 &    0.00000000000014 &    0.00000000000014 &  3.6702E-05 \\
\hline \hline 
 &  0.0 &    0.00167785234899 &    0.00167786399030 &  6.9382E-06  \\
 &   1.0 &    0.00062138979909 &    0.00062139403983 &  6.8246E-06 \\
$ n=300$ &   2.0 &    0.00003157272440 &    0.00003157294168 &  6.8819E-06 \\
 &   3.0 &    0.00000022027431 &    0.00000022027615 &  8.3417E-06 \\
 &   4.0 &    0.00000000021134 &    0.00000000021134 &  8.4873E-06 \\
 &   5.0 &    0.00000000000003 &    0.00000000000003 &  2.6541E-05 \\
\hline 
\end{tabular}
\caption{\small Exact and  approximated values of ${\cL}_{n} u_1(x_1,0,\ldots,0)$ 
and the relative error using $\cL_{4,0.05}^{(n)}$}\label{tabel1}  
\end{center}
\end{small}
\end{table} 

The approximations have been computed on a uniform grid
with step size $h=0.05$ 
using the basis functions
\begin{align*}
\widetilde \eta_{8} (\bdx)=(\pi \cD)^{-n/2} 
\e^{- |\bdx|^2/ \cD} \prod_{j=1}^n  L_{3}^{(1/2)}(x_j^2/ \cD)
\quad \mbox{with} \quad \cD=3.5 \, .
 \end{align*}

In this case one has to determine
\begin{align*}
I_3\Big(\frac{\bdm}{\sqrt{\cD}}\Big)&=\frac{1}{4  } \itg_0^\infty 
\Big(\prod_{j=1}^n  \sum_{k=0}^{3}
        \frac{\e^{\, -m_j^2/(\cD (1+t))}}   {(1+t)^{k}} \, 
L_{k}^{(-1/2)}\Big(\frac{m^2_j}{\cD (1+t)}\Big)
\Big) \, \frac{ dt } {(1+t)^{n/2}} 
\end{align*}
for the integer vectors $\bdm=(m_1,\ldots,m_n) \in [-240,240]^n$,
which reduces by a quadrature
\begin{align*}
I_3\Big(\frac{\bdm}{\sqrt{\cD}}\Big)\approx
\sum_{\ell=1}^R \frac{\omega_\ell}{(1+ \tau_\ell)^{n/2}}\prod_{j=1}^{n}
\sum_{k=0}^{3}
        \frac{\e^{\, -m_j^2/(\cD (1+\tau_\ell))}}   {(1+\tau_\ell)^{k}} \, 
L_{k}^{(-1/2)}\Big(\frac{m^2_j}{\cD (1+\tau_\ell)}\Big)
\end{align*}
to the computation of one vector of length $481$ with the components
\begin{align} \label{1dquad}
\sum_{k=0}^{M-1}\frac{\e^{\, -m^2/(\cD (1+\tau_\ell))}}   
{(1+\tau_\ell)^{k}} \, 
L_{k}^{(-1/2)}\Big(\frac{m^2}{\cD (1+\tau_\ell)}\Big) \,,\quad
m=-240, \ldots,240 \,,
\end{align}
for any quadrature node $\tau_\ell$ and $M=4$.
Let us mention, that for all calculations, reported here and in the
following 
Table  \ref{tabel2},
the same quadrature rule is used 
for computing the integrals $I_3(\bdm/\sqrt{\cD})$.

Table \ref{tabel2} shows
that the cubature method is effective also for much higher
space dimensions. 
We compare the exact values of ${\cL}_n u_2$ and the approximate values  $\cL_{4,h}^{(n)}u_2$
at the same grid points $(x_1,0,\ldots,0) \in \R^n$
for space dimensions $n=10^5,10^6, 2\cdot 10^6$.
The parameters of $\cL_{4,h}^{(n)}$
are $\cD=3.5$ and $h=0.025$.
\begin{table}[!h]
\begin{small}
\begin{center}
\begin{tabular}{r|c|cc|cc|cc} \hline
$n$ &  & $10\,000$ & & $100\,000$ &  & $200\,000$ & \\ \hline
 $x_1$ & exact & abs. error &rel. error & abs. error & rel. error & abs. error & rel. error\\ \hline
 0.0 & -1.0000E-00 & 5.876E-05 &  5.88E-05 &  2.041E-03 &  2.04E-03 &  2.153E-03 &  2.15E-03 \\
 1.0 & -3.6788E-01 & 2.160E-05 &  5.87E-05 &  7.509E-04 &  2.04E-03 &  7.920E-04 &  2.15E-03 \\
 2.0 & -1.8316E-02 & 1.077E-06 &  5.88E-05 &  3.739E-05 &  2.04E-03 &  3.944E-05 &  2.15E-03 \\
 3.0 & -1.2341E-04 & 7.345E-09 &  5.95E-05 &  2.522E-07 &  2.04E-03 &  2.659E-07 &  2.15E-03 \\
 4.0 & -1.1254E-05 & 6.957E-12 &  6.18E-05 &  2.306E-10 &  2.05E-03 &  2.429E-10 &  2.16E-03 \\
 5.0 & -1.3888E-11 & 9.304E-16 &  6.70E-05 &  2.863E-14 &  2.06E-03 &  3.008E-14 &  2.17E-03 \\ 
\hline 
\end{tabular}
\end{center}
\caption{\small Exact values of ${\cL}_{n} u_2(x_1,0,\ldots,0)$, 
absolute and relative errors using $\cL_{4,0.025}^{(n)}$}\label{tabel2}  
\end{small}
\end{table} 
The computing times on a Xeon processor with  $3.4$ GHz
are $0.42$ seconds for  space dimension $n=10\,000$, 
$5.12$ seconds for $n=100\,000$, and $11.63$ seconds for  $n=200\,000$. 

\begin{table}[!h]
\begin{small}
\begin{center}
$M=4$\\[1mm]

\begin{tabular}{r|cc|cc|cc|cc|cc} \hline
 $n $ & $ 3$  & & $ 10$ & &$ 500$ & &$2\,000$ &  &$30\,000$ & \\[1pt] \hline 
$h^{-1} \hspace{-2mm}$ & error  & rate & error  & rate & error  & rate & error  &
rate & error  & rate \\[1pt] \hline 
 5   & 4.99E-05  &          &  6.33E-04  &         &  3.93E-02  &        &  1.34E-01  &        &  3.67E-01  &        \\       
 10  & 4.73E-07  &   6.72   &  4.16E-06  &   7.25  &  2.62E-04  &   7.22 &  1.05E-03  &   6.99 &  1.55E-02  &   4.57 \\
 20  & 2.32E-09  &   7.67   &  1.88E-08  &   7.79  &  1.17E-06  &   7.81 &  4.69E-06  &   7.81 &  7.04E-05  &   7.78 \\
 40  & 9.64E-12  &   7.91   &  7.64E-11  &   7.94  &  4.75E-09  &   7.94 &  1.91E-08  &   7.94 &  2.86E-07  &   7.94 \\
 80  & 4.99E-14  &   7.60   &  4.02E-13  &   7.57  &  2.50E-11  &   7.57 &  1.00E-10  &   7.57 &  1.51E-09  &   7.57 \\
\hline 
\end{tabular}\\[1mm]

$M=3$\\[1mm]

\begin{tabular}{r|cc|cc|cc|cc|cc} \hline
  5   &  1.45E-04  &         &  4.11E-03  &        &  1.98E-01  &  &   3.51E-01  &                 &  3.68E-01  &         \\  
 10   &  5.05E-06  &   4.84  &  9.35E-05  &   5.46 &  6.23E-03  &   4.99 &  2.44E-02  &   3.85      &  2.37E-01  &   0.64  \\
  20  &  9.76E-08  &   5.69  &  1.62E-06  &   5.85 &  1.08E-04  &   5.85 &  4.34E-04  &   5.81      &  6.46E-03  &   5.19  \\
  40  &  1.61E-09  &   5.92  &  2.60E-08  &   5.96 &  1.73E-06  &   5.96 &  6.95E-06  &   5.96      &  1.04E-04  &   5.95  \\
  80  &  2.55E-11  &   5.98  &  4.09E-10  &   5.99 &  2.72E-08  &   5.99 &  1.09E-07  &   5.99 &  1.64E-06  &   5.99  \\ 
\hline 
\end{tabular}\\[1mm]

$M=2$\\[1mm]

\begin{tabular}{r|cc|cc|cc|cc|cc} \hline
  5   & 1.43E-03  &         &  2.89E-02  &         &  3.66E-01  &        &  3.68E-01  &           & 3.68E-01  &        \\  
  10  & 1.04E-04  &   3.78  &  2.32E-03  &   3.64  &  1.29E-01  &   1.51 &  3.02E-01  &   0.28    & 3.68E-01  &   0.00 \\
  20  & 6.99E-06  &   3.90  &  1.55E-04  &   3.91  &  1.04E-02  &   3.63 &  3.98E-02  &   2.92    & 3.02E-01  &   0.28 \\
  40  & 4.46E-07  &   3.97  &  9.83E-06  &   3.98  &  6.66E-04  &   3.96 &  2.67E-03  &   3.90    & 3.81E-02  &   2.99 \\
  80  & 2.80E-08  &   3.99  &  6.17E-07  &   3.99  &  4.18E-05  &   3.99 &  1.68E-04  & 3.99 & 2.51E-03  &   3.92 \\
\hline 
\end{tabular}\\[1mm]

$M=1$\\[1mm]

\begin{tabular}{r|cc|cc|cc|cc|cc} \hline
  5   & 3.73E-02  &         &  1.93E-01  &        &  3.68E-01  &        &  3.68E-01  &             &  3.68E-01  &         \\ 
  10  & 9.29E-03  &   2.00  &  6.56E-02  &   1.56 &  3.68E-01  &   0.00 &  3.68E-01  &   0.00      &  3.68E-01  &   0.00  \\
  20  & 2.31E-03  &   2.01  &  1.79E-02  &   1.88 &  3.51E-01  &   0.07 &  3.68E-01  &   0.00      &  3.68E-01  &   0.00  \\
  40  & 5.75E-04  &   2.00  &  4.56E-03  &   1.97 &  1.99E-01  &   0.82 &  3.52E-01  &   0.07      &  3.68E-01  &   0.00  \\
  80  & 1.44E-04  &   2.00  &  1.15E-03  &   1.99 &  6.50E-02  &   1.61 &  1.99E-01  &   0.82&  3.68E-01  &   0.00  \\
\hline 
\end{tabular}
\caption{\small Absolute errors and approximation rates
for $\cL_n u_2(1,0,\ldots,0)$ using $\cL_{M,h}^{(n)}$}\label{tabel5} 
\end{center}
\end{small}
\end{table} 

In Table \ref{tabel5}  we report 
on the absolute errors and the approximation rates
for the Newton potential $\cL_n u_2(1,0,\ldots,0)=-0.3678794411714423$
in the space dimensions
$n=3,10,500,2\,000$, and $30\,000$.
The approximate values are computed
by the cubature formulas $\cL_{M,h}^{(n)}$ defined by \eqref{cubfM}
for $M=1,2,3,4$,
having the approximation orders \eqref{harmest},
with  $\cD=5$.
We use 
uniform grids of size $h=0.1 \cdot 2^{1-k}$, $k=0,\ldots,4$,
i.e., on the finest grid one has to compute vectors of length $1921$ with the components
\eqref {1dquad} for $m=-960, \ldots,960$. 

The numerical results show that the
high order cubature formulas 
give essentially better approximations
of the Newton potential than second order formulas.
For very high dimensional cases the second order formula fails,
whereas the $8$-th order formula $\cL_{4,h}^{(n)}$ 
approximates with the predicted approximation rate.

\section{Quadrature rules for integral representations }\label{quad}
Here we describe the quadrature rule used to get  separated representations of 
the volume applied to the basis
functions. We give some numerical results on the 
quadrature error for second and fourth order formulas 
approximating the  Newton potential and the potential
of $-\Delta +a^2$.

\subsection{Quadratures}
It is well known that the classical trapezoidal rule
is exponentially converging for certain classes of integrands, for example
periodic functions or rapidly decaying functions on the real line.
For any sufficiently smooth function, say of the Schwartz class $\cS(\R)$,
Poisson's summation formula yields that 
\begin{align*}
h \sum_{k=-\infty}^ \infty f(kh)   = 
\sum_{k=-\infty}^ \infty\hat f \big( \frac{2\pi k}{h}\big) \, .
 \end{align*}
Here  $\hat f$ is the Fourier transform
\[
\hat f(\lambda) =  \itg_{-\infty}^\infty f(x) \e^{-2\pi i x \lambda} \, dx \, .
\]
Thus, 
\[
\itg_{-\infty}^\infty f(x) \, dx - 
h \sum_{k=-\infty}^ \infty f(kh) =
\sum_{k \ne 0}\hat f \big( \frac{2\pi k}{h}\big) \, ,
\]
which indicates that by choosing special substitutions, which transform
the integrand
to a rapidly decaying function with rapidly
decaying Fourier transform, the
trapezoidal rule of step size $h$ can provide very accurate
approximations of the integral for a relatively small
number of quadrature nodes $\{kh\}$.

Additionally, if the integrand is analytic in a strip 
$D_d= \{z \in \C: |\Imm z| < d\}$ such that
\begin{align*}
N(f,D_d)= \| f(\cdot +id)\|_{L_1(\R)} + \| f(\cdot -id)\|_{L_1(\R)} < \infty \, ,
 \end{align*}
then results from Sinc approximation can be applied to 
derive error estimates for the trapezoidal rule.
It was shown in \cite{GHB,HK}, that for doubly exponentially decaying $f$, i.e.,
\begin{align} \label{doubexp}
|f(x)| \le C \exp\big(\!-a \e^{\, b|x|} \! \big) \quad \mbox{for all} \; x \in \R \quad \mbox{with constants }  \; a,b,C > 0 \, , 
 \end{align}
the truncated rule with  $h= \log(2 \pi a  N/b)/(aN)$ satisfies
\[
\Big|\itg_{-\infty}^\infty f(x) \, dx
- h \sum_{k=-N}^N f(kh) \Big| \le C N(f,D_d) e^{\,- 2 \pi a d N/ \log(2 \pi a  N/b)} \, .
\]

\subsection{Newton potentials}
We  make the substitutions
\begin{align} \label{subwald}
t=\e^{\, \xi}\, , \quad \xi=a(\tau + \e^{\tau}) \quad \mbox{and} \quad 
\tau=b(u -\e^{-u})\, ,
 \end{align}
with certain positive constants $a,b$,   proposed in Waldvogel \cite{Wald}. 
Then the integrals \eqref{newton2} are transformed to integrals
over $\R$ with integrands $f(u,\bdx)$
decaying doubly exponentially in $u$.
After the substitution we have
\begin{align*}
&I_1(\bdx)= ab \itg_{-\infty}^\infty 
\frac{(1+\e^{-u})(1 + \exp (b(u -\e^{-u})))\,\phi(u)}
{(1+ \phi(u))^{n/2}} \, \e^{\, - |\bdx|^2/(1+ \phi(u))}\,du
 \end{align*}
with
$ \phi(u)=\exp(ab(u -\e^{-u}) + a\exp(b(u -\e^{-u})))$.
Similarly
\begin{align*}
&I_M(\bdx)= ab \itg_{-\infty}^\infty
\frac{(1+\e^{-u})(1 + \exp (b(u -\e^{-u})))\,\phi(u)}
{(1+ \phi(u))^{n/2}} \,\prod_{j=1}^n g_M(u,x_j) \, du 
\intertext{with the function}
&g_M(u,x)= \e^{\, - x^2/(1+ \phi(u))}
\sum_{k=0}^{M-1} 
\frac{L_{k}^{(-1/2)}\big(x^2/(1+\phi(u) )\big) }
{(1+ \phi(u))^{k}} \, .
 \end{align*}
For the practical application 
it is important
that the   quadrature rule
\begin{equation}\label{trapez}
\itg_{-\infty}^{\infty} f(u,\mathbf{x})\, du
\approx h\, \sum_{k=N_0}^{N_1} f(h k,\mathbf{x})
\end{equation}
approximates $I_M(\bdx)$ 
with prescribed error 
uniformly for $|\mathbf{x}|\leq K$
with a minimal number of summands. 
This number strongly depends on $K$ and on the 
parameters $a,b$ in the transformation \eqref{subwald}.
For the computations in Section~\ref{numres}, which
cover a wide range of dimensions and different orders of the cubature formulas,
we did not try to find optimal parameter sets.
The numbers given in Tables~\ref{tabel1} and \ref{tabel2} were obtained
with
the parameters $a=b=2$ in \eqref{subwald} and $h=.02$, $N_0=-35$,
$N_1=80$ in \eqref{trapez}.
For the results of Table~\ref{tabel5}  
we have chosen $a=6$, $b=5$  and $h=.003$, $N_0=39$, $N_1=250$.
We report in the following on some test to determine optimal parameters for second and fourth
order formulas in low dimensional cases (see also \cite{LMS}).

 \subsubsection{Approximation to the integral $I_1(\mathbf{x})$}

   \begin{figure}[h]
    \centering
       \includegraphics[width=1.9in,height=1.2in]{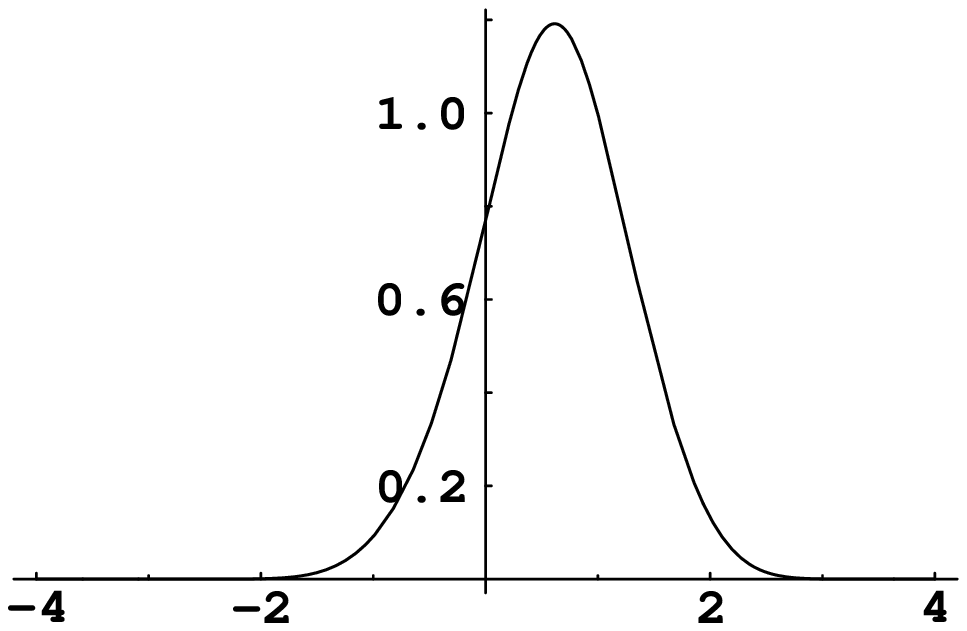} 
        \includegraphics[width=1.9in,height=1.2in]{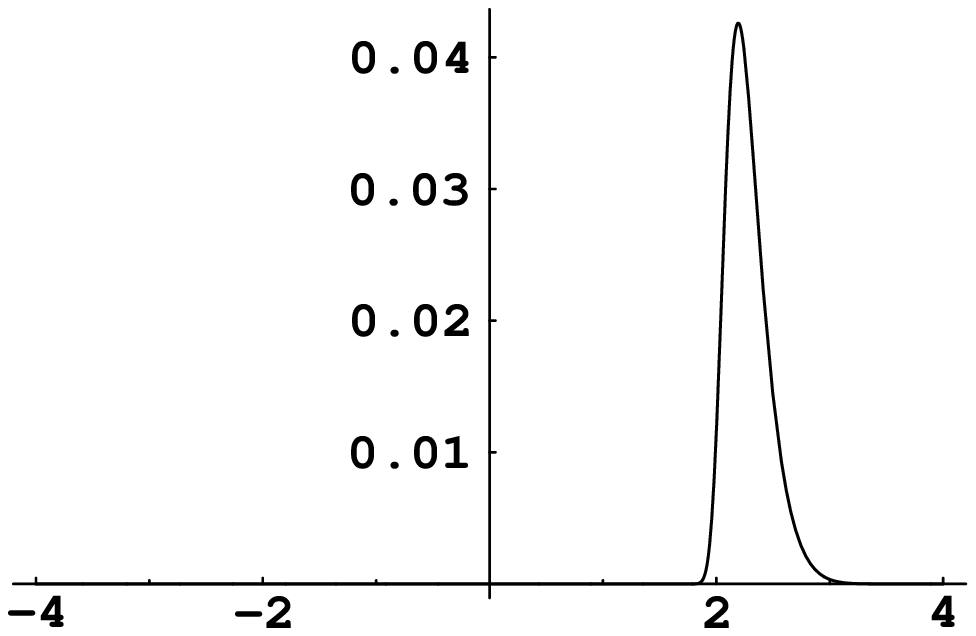} 
       \includegraphics[width=1.9in,height=1.2in]{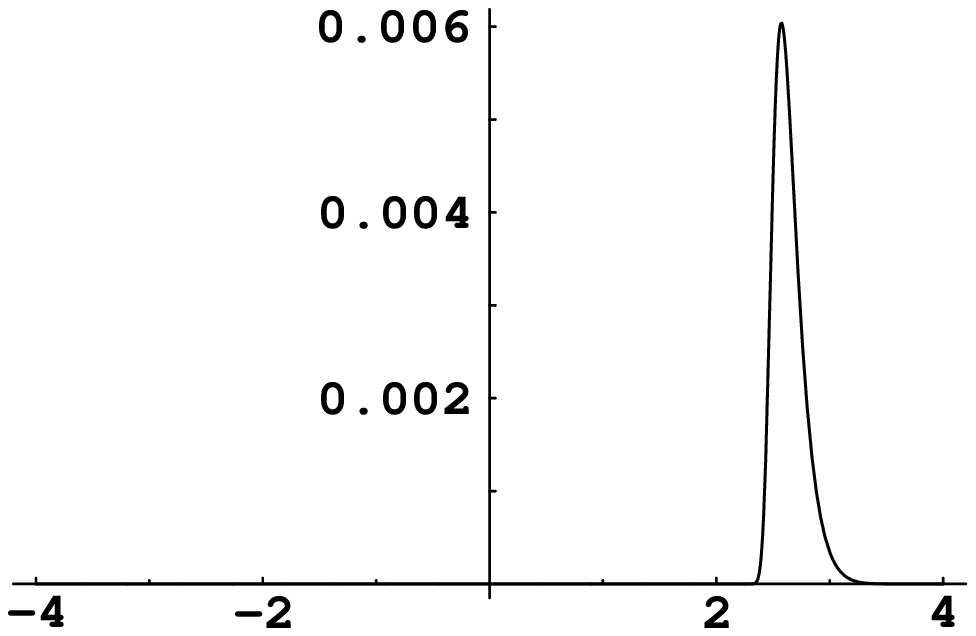}  
        \caption{The plot of the integrand function $f(u,\mathbf{x})$ ($a=b=1$) in 
   $I_1(\mathbf{x})$ for $|\mathbf{x}|=0, 100, 1000$
   (from  the left to the right) in the interval $u\in (-4,4)$.}   
    \label{n3_ab1bis}
    \end{figure}

We assume  in  \eqref{subwald} $a=b=1$.  
Figure~\ref{n3_ab1bis}  illustrates the 
graph of  the integrand function $f(u,\mathbf{x}),\, u \in (-4,4)$, $n=3$,
for different values of $|\mathbf{x}|\leq 10^3$. 
A similar  behavior  holds for other dimensions $n$.
Table \ref{table1} presents  the maximum step $h_0$ 
and the minimum number of quadrature points required to achieve  the relative error $\epsilon$,
uniformly in  $|\mathbf{x}|\in [0,10^3]$  for space dimensions  
$n=3,4,5,6$.  

 \begin{table}[!h]
\begin{small}
\begin{center}
\begin{tabular}{c|cc|cc|cc|cc}
\hline
 & $n=3$ & &$n=4$ & &$n=5$ &&$n=6$ & \\ \hline 
 rel. error & $h_0$& nodes & $h_0$& nodes& $h_0$& nodes& $h_0$& nodes\\
\hline
1E-05  & $0.072$ & $82$&  $0.072$ & $77$ & $0.065 $ & $83$ & $0.059 $ & $90 $\\
1E-07 & $0.055$ & $116$& $0.051$ & $121$& $0.060$  & $96$ &$0.044 $ & $130$\\
1E-09 & $0.043$ & $161$& $0.040$ & $164$& $0.037$ & $169$ & $0.035 $ & $178$\\
1E-11 & $0.036$ & $205$& $0.033$ & $206$& $0.033$ & $200$ & $0.029 $ & $220$\\
 \hline
\end{tabular}
\caption{\small The approximation of $I_1(\mathbf{x})$ for $|\mathbf{x}|\leq 10^3$, with $a=b=1$ 
in \eqref{subwald}}
\label{table1}
\end{center}
\end{small}
\end{table}

\begin{figure}[ht]
   \centering
       \includegraphics[width=1.9in,height=1.2in]{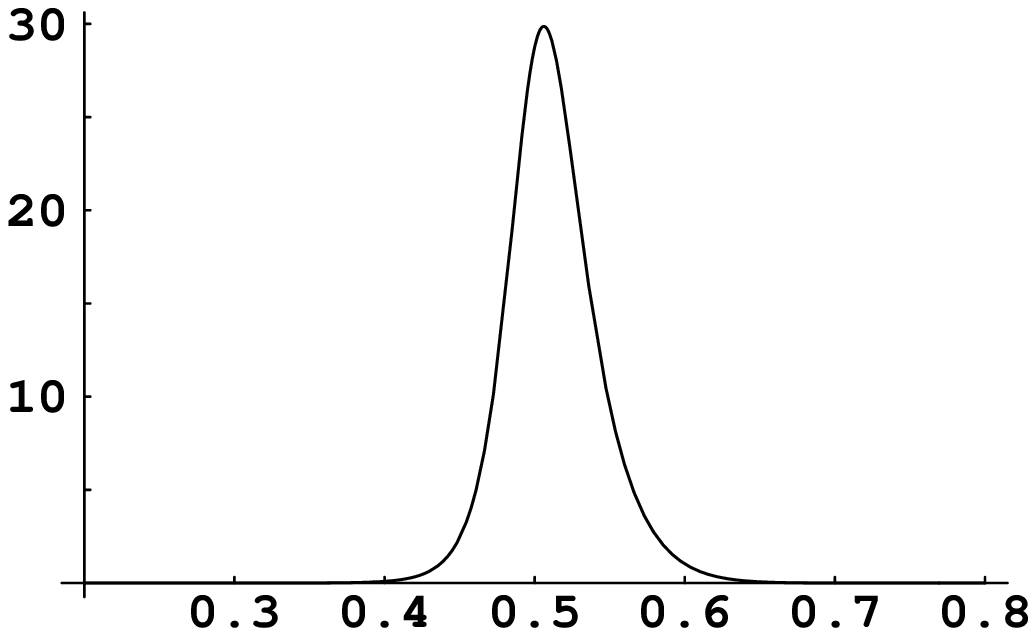} 
        \includegraphics[width=1.9in,height=1.2in]{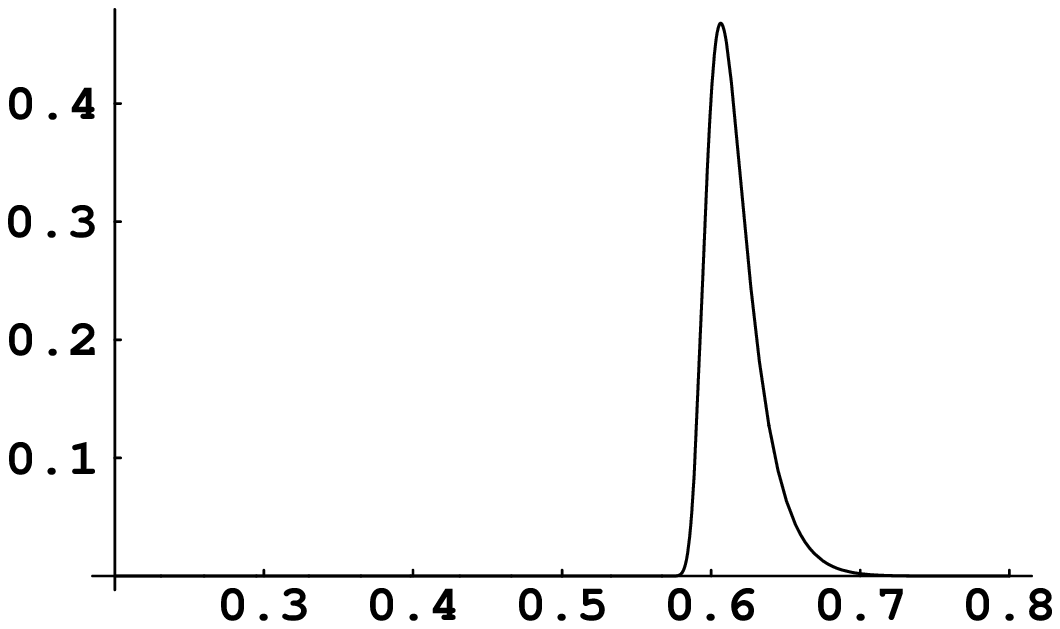} 
       \includegraphics[width=1.9in,height=1.2in]{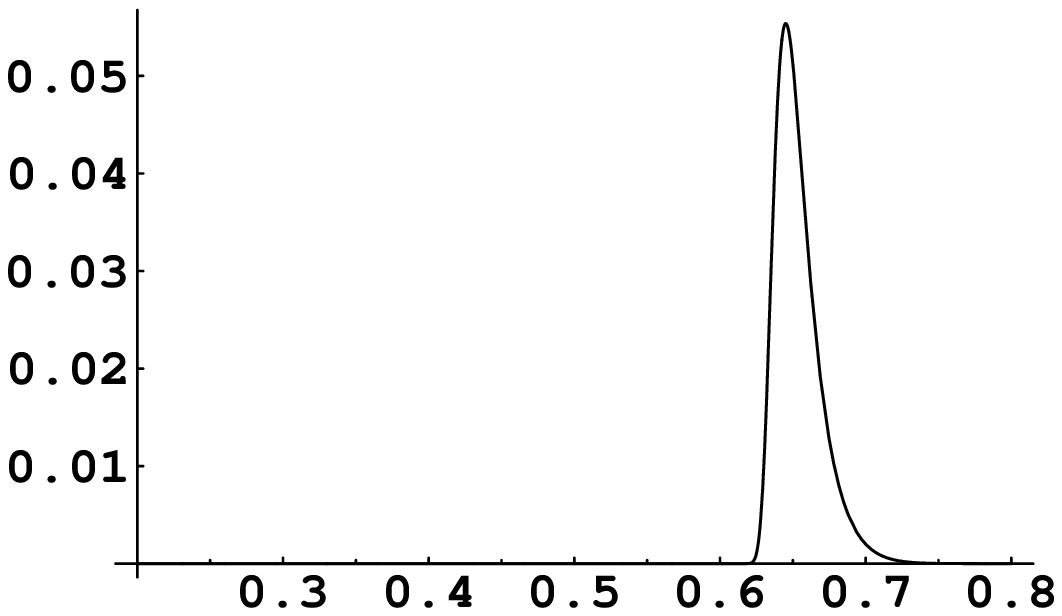}  
        \caption{The plot of the integrand function $f(u,\mathbf{x})$ 
         ($a=6, b=5$) in $I_1(\mathbf{x})$ for 
        $|\mathbf{x}|=0, 100, 1000$
       (from  the left to the right) in the interval $u\in (0,0.85)$.}   
         \label{n3_a6b5bis}
   \end{figure}

It is possible to play with different parameters  $a$ and $b$  in order
to diminish the number of summands  in the quadrature formula. Consider e.g. the case
$a=6$ and $b=5$. Figure \ref{n3_a6b5bis} shows the graph of  $f(u,|\mathbf{x}|)$, $u\in (0,0.85)$
for different values of $|\mathbf{x}|$. 
The numerical results for this quadrature are given in  Table \ref{table2}. 

 \begin{table}[!h]
\begin{small}
\begin{center}
\begin{tabular}{c|cc|cc|cc|cc}
\hline
 & $n=3$ & &$n=4$ & &$n=5$ &&$n=6$ & \\ \hline 
 rel. error & $h_0$& nodes & $h_0$& nodes& $h_0$& nodes &$h_0$& nodes\\
\hline
1E-05 & $0.0077$ & $61$  & $0.0070$ & $83$ & $0.0069$ & $57$   & $0.0058$ & $70 $ \\
1E-07 & $0.0055$ & $111$ & $0.0049$ & $96$ & $0.0046$ & $112$  & $0.0042$ & $117$ \\
1E-09 & $0.0042$ & $170$ & $0.0037$ & $169$& $0.0034$ & $179$  & $0.0037$ & $158$ \\
1E-11 & $0.0034$ & $247$ & $0.0033$ & $200$ & $0.0031$ & $221$  & $0.0028$ & $242$ \\
 \hline
\end{tabular}
\caption{\small The approximation of $I_1(\mathbf{x})$ for $|\mathbf{x}|\leq 10^3$, with  $a=6,\, b=5$ 
in \eqref{subwald}}
\label{table2}
\end{center}
\end{small}
\end{table}

 \subsubsection{Approximation to the integral  $I_2(\mathbf{x})$}

  Next, we discuss the computation of the integral 
 \[
 I_2(\mathbf{x})=\itg_0^\infty \prod_{j=1}^n {\rm e}^{-x_j^2/(t+1)}
 \Big(\frac{3+2t}{2(1+t)} -\frac{x_j^2}{(1+t)^2}\Big)\, 
\frac{dt}{(1+t)^{n/2}} 
 \]
using the variable transformation \eqref{subwald} and the trapezoidal rule \eqref{trapez}, 
for $n=3$ and $n=4$. 
The graphs of the integrands $f(u,\mathbf{x})$ are very similar to the
case of $I_1(\mathbf{x})$.
In the numerical results, for the sake of simplicity,  we assumed 
$\mathbf{x}=(x,x,x)$, with $|\mathbf{x}|\leq 10^3$.
The maximum step $h_0$ 
and the minimum number of quadrature points required to achieve  the 
relative error $\epsilon$ uniformly in  $|\mathbf{x}|\in [0,10^3]$
are presented in Table \ref{table3}
for the cases $a=b=1$ and  $a=6,\, b=5$.

 \begin{table}[!h]
\begin{small}
\begin{center}
\begin{tabular}{c|cc|cc|cc|cc}
\hline
 & $n=3$ & &$n=4$ & &$n=3$ &&$n=4$ & \\ \hline 
 rel. error & $h_0$& nodes & $h_0$& nodes& $h_0$& nodes &$h_0$& nodes\\
\hline
1E-05 & $0.072$ & $82$  & $0.072$ & $77$  & $0.0077$ & $63$  & $0.0074$ & $57$  \\
1E-07 & $0.055$ & $118$ & $0.051$ & $121$ & $0.0052$ & $114$ & $0.0046$ & $120$ \\
1E-09 & $0.043$ & $163$ & $0.040$ & $163$ & $0.0042$ & $175$ & $0.0037$ & $175$ \\
1E-11 & $0.036$ & $204$ & $0.033$ & $206$ & $0.0034$ & $234$ & $0.0033$ & $222$ \\
 \hline
\end{tabular}
\caption{\small The approximation of $I_2(\mathbf{x})$ for $|\mathbf{x}|\leq 10^3$, 
with $a=b=1$ (left) and  $a=6, \,b=5$ (right)}
\label{table3}
\end{center}
\end{small}
\end{table}

\subsection{Potentials of advection-diffusion operators}
\begin{figure}[!h]
       \includegraphics[width=1.9in,height=1.2in]{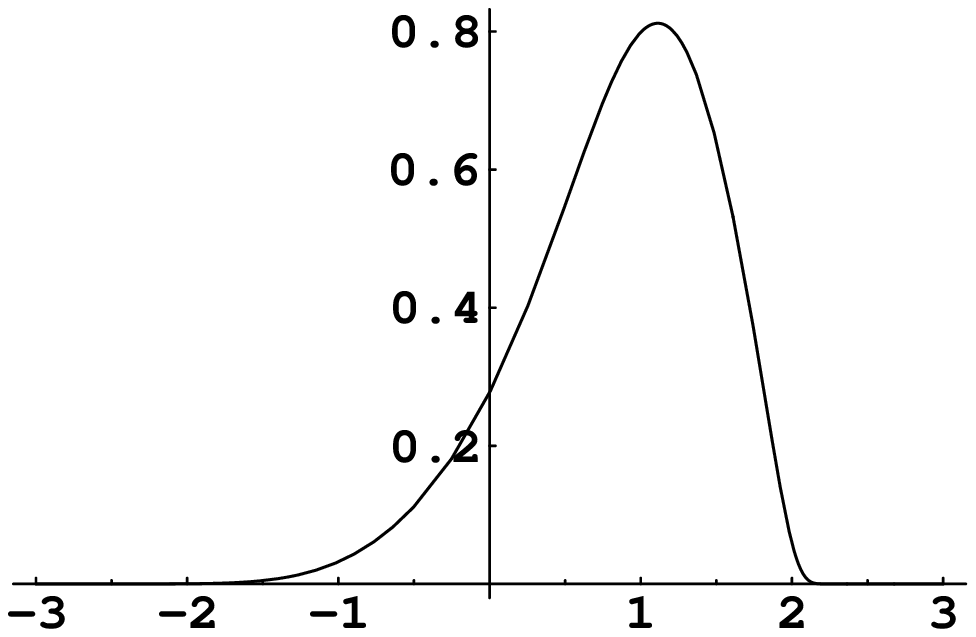} 
        \includegraphics[width=1.9in,height=1.2in]{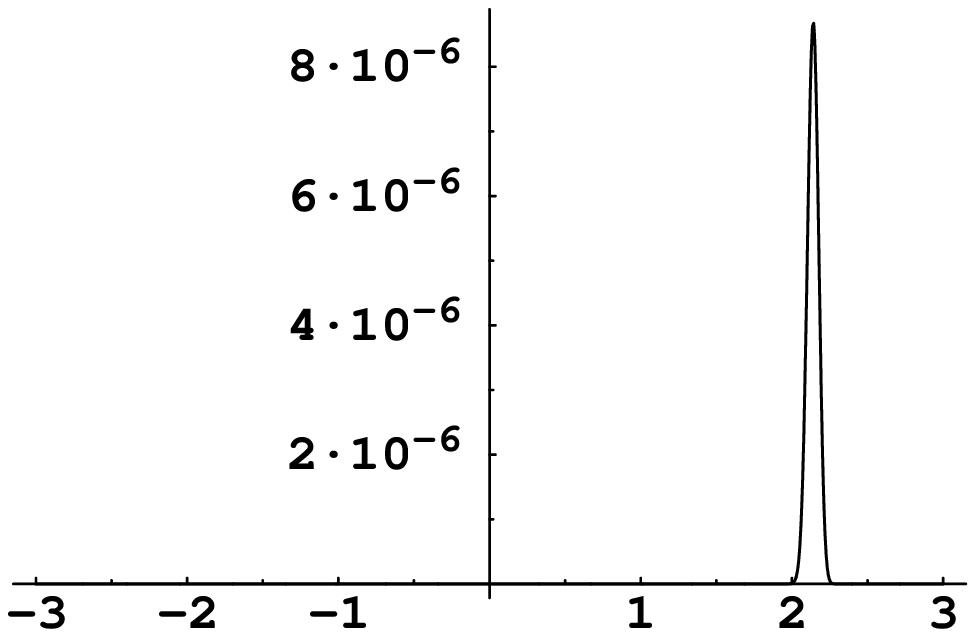} 
       \includegraphics[width=1.9in,height=1.2in]{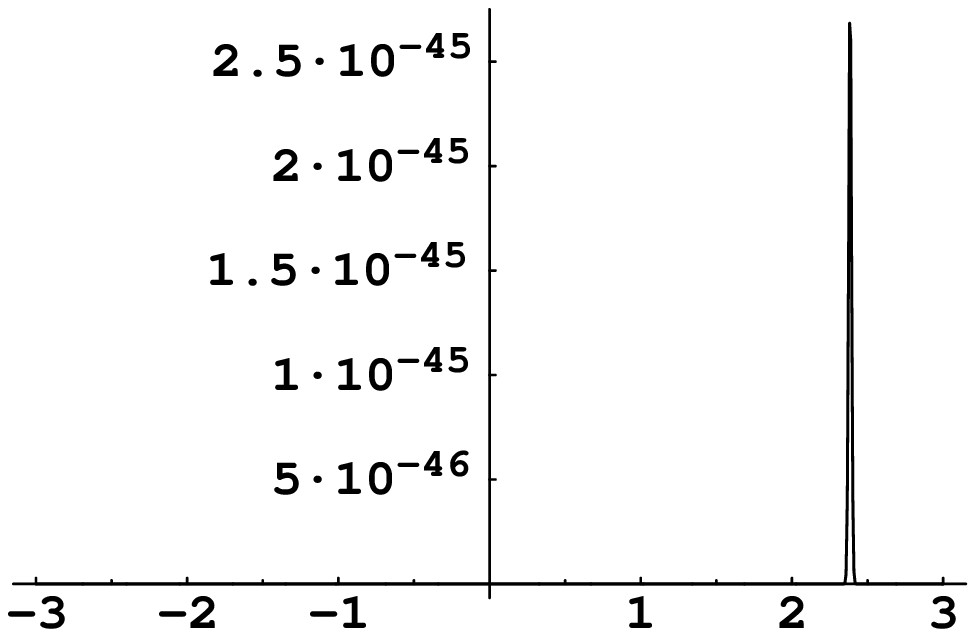}   
        \caption{The plot of the integrand function 
        $f(u,\mathbf{x})$ in $K_1(\mathbf{x}),\> a^2=0.01, b=1$ 
         for $|\mathbf{x}|=0, 100, 1000$
          (from  the left to the right) in the interval $u\in (-3,3)$.}   
         \label{helm1}
  \vskip5pt

       \includegraphics[width=1.9in,height=1.2in]{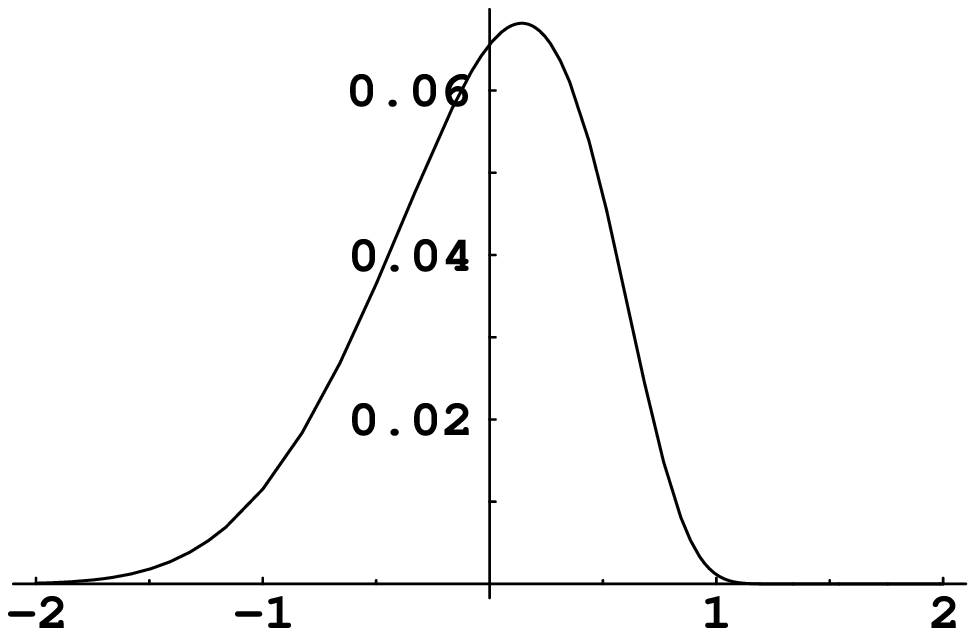} 
        \includegraphics[width=1.9in,height=1.2in]{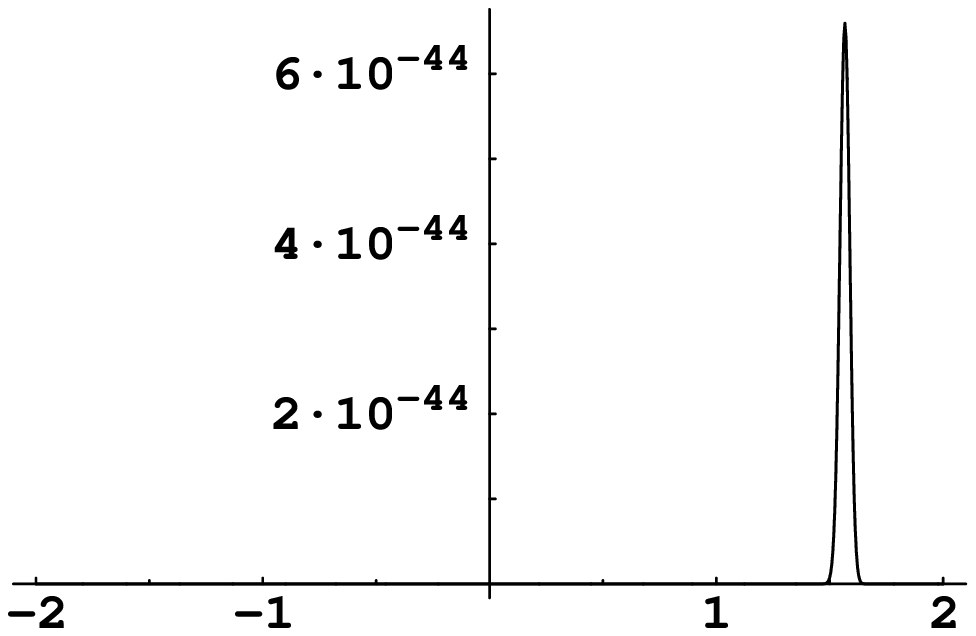}  
           \includegraphics[width=1.9in,height=1.2in]{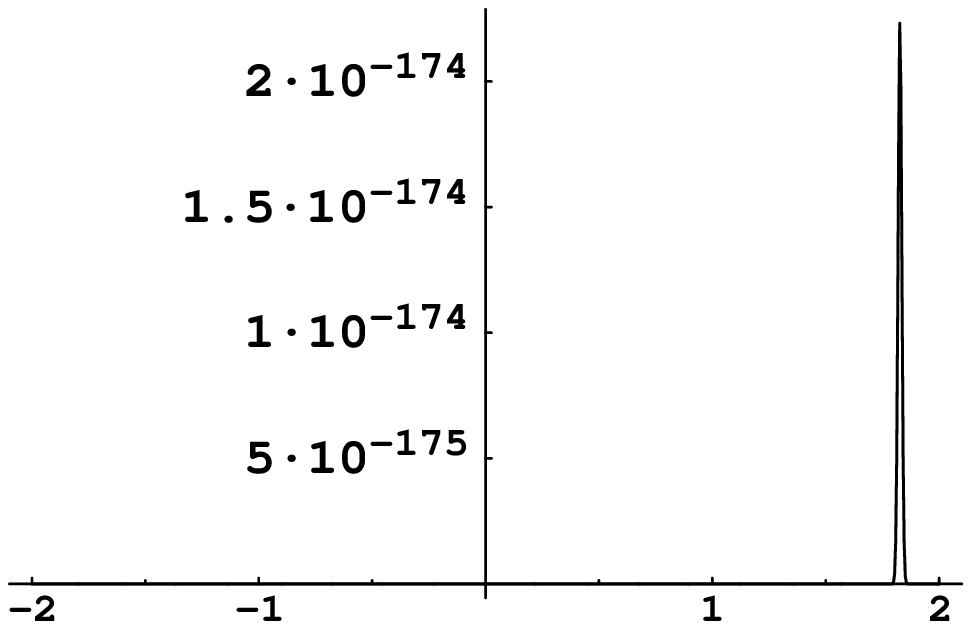}    
        \caption{The plot of the integrand function 
         $f(u,\mathbf{x})$ in $K_1(\mathbf{x}),\> a^2=4, b=1$ 
        for $|\mathbf{x}|=0, 50, 200$
        (from  the left to the right) in the interval $u\in (-2,2)$.}   
         \label{helm2}
    \end{figure}

We consider the special case $\mathbf{b}=0$, $c=a^2>0$, which leads to 
the one-dimensional representation of the volume potential applied to $\prod \widetilde \eta_{2M}$
\begin{align*}
K_M(\bdx)=\itg_0^\infty\e^{-a^2  t/4} 
\prod_{j=1}^n \Big( \sum_{k=0}^{M-1}
        \frac{\e^{\, -x_j^2/(1+t)}}   {(1+t)^{k}} \, 
L_{k}^{(-1/2)}\Big(\frac{x^2_j}{1+t}
\Big) \Big)\, \frac{dt}{(1+t)^{n/2}} \, , \quad M \ge 1
 \end{align*}
(cf. \eqref{intadv}, \eqref{intadvM}). 
To get doubly exponential integrands we make the substitution
\begin{align} \label{subadv}
t=\phi(u) \quad \mbox{with }  \phi(u)=\exp(b(u-\exp(-u)) \, , \; b>0 \, ,
 \end{align}
and apply the trapezoidal rule to
\begin{align*}
b\itg_{-\infty}^\infty
\frac{ \e^{-a^2  \phi(u)/4} (1+\e^{-u})\,\phi(u)}{(1+\phi(u))^{n/2}}  \e^{-|\bdx|^2/(1+ \phi(u))} 
\prod_{j=1}^n
\sum_{k=0}^{M-1}
        \frac{L_{k}^{(-1/2)}(x^2_j/(1+\phi(u)))}   {(1+\phi(u))^{k}} \, du \, .
 \end{align*}
Figures~\ref{helm1} and \ref{helm2} illustrate the 
graph of  the integrand  $f(u,\mathbf{x})$ 
of $K_1(\bdx)$, $n=3$, for different values of $|\mathbf{x}|$
and $a^2=0.01$ and  $a^2=4$.

In Table~\ref{table5} we present the maximum step $h_0$ 
and the minimum number of quadrature points in formula \eqref{trapez} 
for the integral $K_1(\mathbf{x})$  for $n=3$,
in the cases $a^2=0.01,0.1, 1, 4$.  
We have chosen  $b=1$ in the substitution \eqref{subadv}.  

 \begin{table}[!h]
\begin{small}
\begin{center}
\begin{tabular}{c|cc|cc|cc|cc}
\hline
$a^2$& & $0.01$ & &$0.1$ & &$1$ &&$4$  \\ \hline 
 rel. error & $h_0$& nodes & $h_0$& nodes& $h_0$& nodes &$h_0$& nodes\\
\hline
1E-05  & $0.58 $ & $20 $ & $0.58 $ & $17$ & $0.48 $ & $15 $ &  $0.44 $ & $13 $ \\
1E-07  & $0.47 $ & $25 $ & $0.42 $ & $16$ & $0.38 $ & $20 $ &  $0.36 $ & $17 $ \\
1E-09  & $0.39 $ & $32 $ & $0.40 $ & $25$ & $0.37 $ & $22 $ &  $0.31 $ & $21 $ \\
1E-11  & $0.30 $ & $43 $ & $0.29 $ & $36$ & $0.29 $ & $28 $ &  $0.27 $ & $25 $ \\
1E-13  & $0.26 $ & $50 $ & $0.25 $ & $43$ & $0.25 $ & $34 $ &  $0.25 $ & $29 $ \\
 \hline
\end{tabular}
\caption{\small The approximation of $K_1(\mathbf{x})$ with different values
  of $a^2$ for $|\mathbf{x}|\leq 10^3$}
\label{table5}
\end{center}
\end{small}
\end{table}

\subsubsection*{Acknowledgments}
This research was partially supported by the UK and Engineering
and Physical Sciences Research Council via the grant EP/F005563/1.

The authors would like to thank B. Khoromskij
for valuable discussions concerning fast computations of 
high-dimensional problems.


\begin{thebibliography}{100}
\bibitem{Abr}
{M.~Abramowitz and I.~A. Stegun}, {\em Handbook of Mathematical Functions},
  Dover Publ., New York, 1968.

\bibitem{BCFH}
G.~Beylkin, R.~Cramer, G.~Fann, and R.~J.~Harrison,
Multiresolution separated representations of singular and weakly
singular operators,
Appl. Comp. Harm. Anal. 23 (2007), 235--253.

\bibitem{BM1}
G.~Beylkin and M.~J.~Mohlenkamp,
Numerical operator calculus in higher dimensions,
Proc. Nat. Acad. Sci. USA 99 (2002), 10246--10251.

\bibitem{BM2}
G.~Beylkin and M.~J.~Mohlenkamp,
Algorithms for numerical analysis in high dimensions,
SIAM J. on Scient. Comp. 26, 6 (2005), 2133--2159. 

\bibitem{BHK} 
C. Bertoglio, W.~Hackbusch, and B.~N.~Khoromskij,
Low rank tensor-product approximation of projected Green 
kernels via sinc-quadratures,
Preprint 79, MPI MIS 2008.

\bibitem{GHB}
I.~P.~Gavrilyuk,  W.~Hackbusch, and B.~N.~Khoromskij,
Hierarchical tensor-product approximation of the inverse and related
operators in high-dimensional elliptic problems, Computing 74 (2005), 131--157.

\bibitem{H} W.~Hackbusch, Efficient convolution with the Newton
  potential in {\em d} dimensions. Numer. Math. 110 (2008), 449--489.

 
\bibitem{HK} W.~Hackbusch and B.~N.~Khoromskij,
Tensor-product approximation to 
multi-dimensional integral operators and Green's functions,
SIAM J. Matrix Anal. Appl.  30, 3 (2008),  1233 - 1253.

\bibitem{Khor} 
B.~N.~Khoromskij,
Fast tensor approximation of multi-dimensional convolution with linear scaling,
Preprint 36, MPI MIS 2008.

\bibitem{LMS} F.~Lanzara, V. Maz'ya, and G. Schmidt, 
Tensor product approximations of high dimensional potentials,
Preprint 1403, WIAS Berlin 2009.

\bibitem{Maz}
 V. Maz'ya,
  {Approximate approximations}, in {\em The Mathematics of Finite
  Elements and Applications. Highlights 1993}, J.~R. Whiteman, ed., Wiley \&
  Sons, Chichester, 1994, 77--104.

\bibitem{MS2}
{V. Maz'ya and G. Schmidt}
   { ``Approximate
   Approximations'' and the cubature of potentials}, 
Rend. Mat. Acc. Lincei,
    6 (1995), s. 9,   161--184.

\bibitem{MSbook}
 {V. Maz'ya and G. Schmidt}, {\em Approximate
Approximations},  Math. Surveys and Monographs vol. 141, AMS 2007.

\bibitem{MS5} {V. Maz'ya and G. Schmidt},
 { Potentials of Gaussians and approximate wavelets},
Math. Nachr.  280  (2007),  no. 9-10, 1176--1189.

\bibitem{Wald}
J. Waldvogel, 
Towards a general error theory of the trapezoidal rule.
Approximation and Computation 2008.
Nis, Serbia, August 25-29, 2008:\\
\url{http://www.math.ethz.ch/~waldvoge/Projects/integrals.html}

\end{thebibliography}
\end{document}